\documentclass[fleqn,10pt]{wlscirep}
\usepackage[utf8]{inputenc}
\usepackage[T1]{fontenc}

\graphicspath{{./}{./images/}}
\graphicspath{{./}{./figures/}}
\usepackage{multirow}
\newcommand*\diff{\mathop{}\!\mathrm{d}}

\title{Search strategy in a complex and dynamic environment: the MH370 case}

\author[1,*]{Stefan~Ivić}
\author[2]{Bojan~Crnković}
\author[3]{Hassan~Arbabi}
\author[4]{Sophie~Loire}
\author[5]{Patrick~Clary}
\author[6]{Igor~Mezić}
\affil[1]{Faculty of Engineering, University of Rijeka, Rijeka, Croatia}
\affil[2]{Department of Mathematics, University of Rijeka, Rijeka, Croatia}
\affil[3]{Department of Mechanical Engineering, Massachusetts Institute of Technology, Cambridge, USA}
\affil[4]{Bruker Nano Surfaces, Santa Barbara, USA}
\affil[5]{Clary is with Oregon State University, Corvallis, USA}
\affil[6]{Department of Mechanical Engineering and The Center for Control, Dynamical Systems and Computation, University of California, Santa Barbara, USA}

\affil[*]{Corresponding author, stefan.ivic@riteh.hr}

\begin{abstract}
Search and detection of objects on the ocean surface is a challenging task due to the complexity of the drift dynamics and lack of known optimal solutions for the path of the search agents. This challenge was highlighted by the unsuccessful search for Malaysian Flight 370 (MH370) which disappeared on March 8, 2014. In this paper, we propose an improvement of a  search algorithm rooted in the ergodic theory of dynamical systems which can accommodate complex geometries and uncertainties of the drifting search areas on the ocean surface. We illustrate the effectiveness of this algorithm in a computational replication of the conducted search for MH370. In comparison to conventional search methods, the proposed algorithm leads to an order of magnitude improvement in success rate over the time period of the actual search operation. Simulations of the proposed search control also indicate that the initial success rate of finding debris increases in the event of delayed search commencement. This is due to the existence of convergence zones in the search area which leads to local aggregation of debris in those zones and hence reduction of the effective size of the area to be searched. 
\end{abstract}
\begin{document}

\flushbottom
\maketitle
%
%
\thispagestyle{empty}


\section{Introduction}

The problem of search and detection is a crucial part of aerial and maritime operations. The importance of this task was accentuated by the tragic mystery of the missing Malaysian Airline Flight, known as MH370 \cite{mcnutt2014hunt}. The surface search operation that began in the aftermath of the disappearance involved as many as 40 ships and aircraft sweeping an area of over 4,600,000 square kilometers in the Indian Ocean but ended without finding any conclusive traces of the passengers or the plane wreckage \cite{bureau2017operational}. 

There are two significant sources of difficulty in designing search operations on the surface of the ocean. First is the high level of uncertainty associated with the assessment of the search target area. In the case of MH370, the target areas were continuously reassigned due to a stream of incoming information from satellite data analysis, suspicious debris sightings, and detection of signals suspected to be from the underwater locator beacon. This type of uncertainty is further compounded by the complex drift dynamics on the surface of the ocean. The unsteady nature of ocean currents amplifies the errors in the estimation of the splash location in the days lapsed between the splash and the start of the search operation. In addition, given that drift models also have inaccuracies, the final computational estimates for the search targets typically have a considerable amount of uncertainty. A more detailed studies are published on topics of complexity and uncertainties of MH370 path and splash estimation \cite{ashton2015search,davey2016bayesian}, and debris drift affected by the stirring of the ocean \cite{garcia2015dynamical,jansen2016drift,trinanes2016analysis,nesterov2018consideration}.

The second source of difficulty is that even when the target area for search is known and the associated uncertainty is well-characterized, e.g. in the form of a probability distribution, the problem of the design of the search agents’ paths generally does not have a known optimal solution. The conventional approach is to divide the target area into subdomains and assign them to different agents which will sweep each subdomain in organized paths, e.g., parallel lines. This so-called “lawnmower” strategy is easy to plan and implement for simple geometries, but it cannot be readily extended to complex geometries that arise from the surface drifts or areas with non-uniform distribution of likelihood for finding targets. The lawnmower technique also lacks the flexibility required for real-world operations: for example, it would require reassignment of search areas and agents in case of an agent going astray.

We propose a multi-agent motion control method called modified Dynamic Spectral Multi-scale Coverage (mDSMC), for search and detection of objects in dynamically evolving environments such as the ocean surface. This algorithm combines the classical theory of optimal search \cite{koopman1946search,stone1976theory} with concepts from ergodic theory and is capable of accommodating complex geometries, non-uniform distributions and the instantaneous drift of targets during the search. The algorithm is particularly designed to downplay the role of small spatial characteristics of the target area which leads to a significant increase in the success rate of the search operation. To show the promise of this algorithm for real-world applications, we compare its performance with conventional search techniques in a computational replication of the surface search for MH370.

\section{Related research}

The design of search agents' path is a multi-agent control problem whose objective function involves the estimation of a time evolving probability distribution.
Although considering the dynamics is a key feature of successful oceanic search, there are a lot of relevant and interesting details that can be found in papers discussing both static and dynamic search problems. The following overview of related literature presents the relevant publications grouped to approaches considering static targets and the ones where targets are dynamic.

\subsection{Stationary targets}

In \cite{yao2017optimal} a route planning for Unmanned Aerial Vehicle (UAV) in stationary target search mission over a river region is considered. Sub-regions along the river are extracted using a Gaussian mixture model (GMM) and prioritized heuristically with an approximation insertion (AI) approach. The optimal routes are obtained maximizing stationary target detection in a given time window. This approach is extended to the search with multiple UAVs in \cite{yao2017gaussian} where Receding Horizon Control (RHC, also known as Model Predictive Control or MPC) is employed for finding the optimal paths.
In \cite{miller2015ergodic} RHC is employed to solve ergodic exploration of distributed information.
It is shown that optimization-based approach is suitable for both local and global search density i.e. whether the information is localized or diffused. 
The optimization-based path-planner, presented by Gramajo and Shankar \cite{gramajo2017efficient}, optimizes the search for given energy consumption and maneuverability constraints.

An advanced lawnmower-type control algorithm presented in \cite{kingston2016automated} considers the dynamics of the UAV, the spatial scope and the quality of the sensor.
Another lawnmower-based search improvement \cite{lin2014hierarchical} finds paths that approximate the payoff of an optimal solution based on partial detection in the form of a task difficulty map. The algorithm uses the mode goodness ratio heuristic that uses a Gaussian mixture model to prioritize search subregions.

A reconnaissance mission control presented by Wang et al. \cite{wang2019reconnaissance} focuses on particle swarm optimization for UAV swarm path planning. The optimization considers different search strategies and optimization objectives, and constraints regarding the mobility and communication of the UAV swarm.
A machine learning approach for the search and rescue in indoor environments using UAVs has been investigated in \cite{kulkarni2020uav}. The reinforced learning is used to locate the trapped victim by sensing the RF signals emitted from the victim's smart device.
A heat equation driven area coverage (HEDAC) control method \cite{ivic2017ergodicity} has been employed for heterogenous multi-agent search in uncertainty conditions \cite{ivic2019motion}. An exact probabilistic model and state-of-the-art control method ensures a near-optimal performance for stationary targets.  
An algorithm called layered search and rescue (LSAR) is employed in \cite{alotaibi2019lsar} for multi-agent search and rescue missions in order to minimize the search time while finding the maximum number of victims.


\subsection{Non-stationary targets}

An Ant Colony Optimization (ACO) is used in \cite{perez2018ant} to determine paths of multiple UAVs that ensure minimal time needed to find moving targets.
A specialized multi-UAV sea area search map is presented in \cite{yue2019new}, where the target probability map (TPM) was designed to handle uncertainties caused by dynamic targets. The TPM is used as a pheromone map in an improved multi-ant colony algorithm - a derivation of ACO algorithm. 
The main drawback of this heuristic approach is that the resulting paths are straight-segmented due to the discrete nature of the method.
Another heuristic approach is presented in \cite{du2019evolutionary} where several algorithms are combined in a multi-UAV search procedure for missing persons with the time-varying distribution of target location probabilities. 
The algorithm used for evolving a population of main solutions is a hybrid evolutionary optimization method, while each UAV path is finely optimized using Tabu Search method.

A distributed cooperative multi-UAV search based on the prediction of non-stationary targets is proposed in \cite{ru2015distributed}.
Target existence probability is updated using detection information of on-board sensors while the UAV's are directed using RHC for maximizing search performance.
Ergodic exploration using RHC, presented in \cite{mavrommati2017real}, plans a real-time motion that locally optimize ergodicity with respect to the dynamic information density. 
RHC is also utilized in \cite{hu2017optimal} where targets are partially aware of the UAVs locations and are trying to escape them. Although this is not the case in MH370 search, it provides an interesting game-theoretic view of the search problem \cite{li2017potential}.   

A mission planner for locating and tracking harmful ocean debris with UAVs is presented in \cite{rubio2004adaptive}.
Actual weather data and predicted icing conditions and their impact on the UAV performance are taken into account in the search simulations. Market-based cooperation strategies for multiple UAVs and evolutionary computation techniques are utilized for conducting low-cost oceanic search missions.

A path planning algorithm \cite{sujit2009coordination}  governs the coordinated search using unmanned aerial vehicles and autonomous underwater vehicles. The proposed control framework is tested on lawnmower-style multi-agent search simulations. 
A search control for autonomous underwater vehicles is studied in \cite{huang2014dynamic} where non-stationary targets driven by ocean current are considered. A search for multiple target locations in an unsteady ocean environment is achieved with a self-organizing map (SOM) neural network. Considering target movement due to ocean currents, an optimal search agent path is found using a velocity synthesis approach.

Yau and Chung \cite{yau2012search} investigate the application of linear search and discrete myopic search, coupled with surrogate ocean models, for locating a drifting object.
A motion and camera control algorithm for multiple UAVs proposed by Perez et al. \cite{perez2019minimum} minimizes the expected detection time of a nondeterministically moving target of an uncertain initial position. The method is tested on 3 real-world inspired scenarios including a drifting boat by the coast.
Multi-UAV search for moving targets in an unknown environment is established in \cite{yue2019novel} using a reinforced learning scheme. The technique is applied in a search of moving ships at the sea.



\section{Search in dynamic environments}
We will consider a rectangular domain $S\subset \mathbb{R}^2$ with boundary $\partial S$ which is large enough to contain all potential debris floating in the sea and agent trajectories for the entire search period. The algorithm that we use here is a new variant of the Spectral Multi-scale Coverage (SMC) algorithm which was proposed in \cite{mathew2011metrics}. This algorithm was extended to dynamic environments in \cite{mathew2010uniform} (called DSMC) and its practical feasibility was shown in \cite{mathew2013experimental}. We give a short introduction to DSMC algorithm and emphasize the alterations and that leads to the algorithm introduced in this paper called modified DSMC (mDSMC).

From a mathematical viewpoint, a search problem is formulated around two central concepts. First is the probability distribution function of the search targets (i.e. survivors, debris, etc.) which denotes the likelihood of finding targets in each subset of the domain.  This distribution, denoted here by $p^t$, is dynamic and evolves in time according to the drift equation in an incompressible flow:
\begin{equation}
\frac{\partial p^t}{\partial t} + v \nabla p^t = 0, \qquad \forall x\in S,
\label{eq:drift}
\end{equation}
where $v$ is the ocean surface velocity and the initial distribution  $p^{t_0}$ specifies the estimated splash location and the uncertainties therein at the estimated time of the splash $t_0$.  

We denote a flow map $\mathcal{T}^{t_1,t_2}: \mathbb{R}^2 \mapsto \mathbb{R}^2$, which maps a location of a sample (or a target) at time ${t_1}$ to its location at time ${t_2}$ using the unsteady velocity field $v$.
The flow map $\mathcal{T}$ describes the evolution of target and sample positions on the ocean surface which is a basis for the use of Lagrangian methods. The advection of samples placed on the initial splash location through flow map composition yields numerical solutions almost devoid of artificial diffusion and allows for direct parallelization.

The second concept is the coverage distribution $c^t$, which reflects the spatial distribution of the conducted search time up to time $t$. In other words, regions with higher coverage have been traversed by search agents more often. Assuming $N$ search agents are dispatched at time 0 for a ceaseless search operation, and the path of each agent is denoted by $z^t_i$ with $i=1,2, \ldots, N$, the coverage distribution is given by 
\begin{equation}
c^t  =\sum_{i=1}^N \int_0^t  \delta({x} - \mathcal{T}^{\tau,t}({z}_i^\tau)) \diff \tau ,
\label{eq:coverage}
\end{equation} 
where $\delta$ is the Dirac delta function. The coverage can be thought of as a distribution
with support on the points visited by the agents, but
evolved forward in time by drift dynamics of the ocean surface.
Note that the total search coverage at time $t$ can be computed by integrating the above density over the whole area of search $S$,
\begin{equation}
\int_{S}c^t \diff x = N t.
\label{eq:total_coverage}
\end{equation}
When considering the concept of coverage, it is useful to model the search process by taking into account the properties of the detection system used in the search. The agent is able to detect targets near the current location but the detection ability decreases as the distance between the agent and the target increases. A simple model of this detection process and data acquisition can be described using a positive smooth compact support Radial Basis Functions (RBF) $\phi$ and $\phi_\sigma$ that satisfy:
\begin{equation*}
\phi_\sigma({x})=\sigma^{-2}\phi\left(\frac{{x}}{\sigma}\right)\;
\textmd{ and }
\int_{\mathbb{R}^2} \phi({x})\diff x = 1,
\end{equation*}
where $\sigma$ is a positive scaling factor that controls the visual detection range. The RBF models the measurement density of the visual detection system. To incorporate the RBF model we must also modify the coverage and define a smooth modified coverage function which represents the approximation of total measurement density in the domain up to time $t$:
\begin{equation}
c_\sigma^t(x) =\phi_\sigma(x)*c^t,
\label{eq:c_smoothing}
\end{equation}
which also satisfies the condition:
\begin{equation}
\int_S c^t_\sigma \diff x = N t.
\label{eq:total_smooth_coverage}
\end{equation}
Furthermore, one can observe that these two definitions are consistent by passing to the limit:
\begin{equation*}
\lim_{\sigma\to 0^+}c^t_\sigma(x) = c^t.
\end{equation*}
If the coverage \eqref{eq:c_smoothing} is integrated over $S$ and the support of the RBF function is outside $S$, one must rescale the coverage to keep growth of total covearge linear in time.

This modification of coverage is not an exact model of coverage drift since it does not allow for the stretch and contraction of coverage support during the drift (i.e. $supp(c_\sigma^t)\subset supp(c^t)+supp(\phi_\sigma)$, where $+$ denotes the Minkowski addition) but it models the search process more faithfully compared to \eqref{eq:coverage}. Furthermore, it makes implementation straightforward and allows for real-time computation.

The goal of the search operation planning is to design the paths of the search agents such that the probability of detecting targets is maximized over a given time window. From \cite{koopman1957theory}, the probability of the target being detected by time $t$ is 
\begin{equation}
P_d = \int_{S} p^t \mathcal{F}(c^t_\sigma({x})) \diff x
\label{eq:detection_prob}
\end{equation}
where $\mathcal{F}$ is called the detection function. To be more precise, $\mathcal{F}(c^t_\sigma({x}))$ denotes the (conditional) probability of the target being detected given the target is at $x$. The prevalent choice of detection function, which comes from the assumption of locally random search \cite{koopman1956theory,stone2016search}, is the exponential saturation. Mathematically speaking, this law is stated as
\begin{equation*}
\mathcal{F}(c) = 1 - e^{-c}.
\end{equation*}
This means that if the agent and the target visit the same location at the same time, there is still some probability that the agent misses the target, but this probability decays exponentially with each visit.  

The optimal coverage distribution that maximizes the probability in \eqref{eq:detection_prob}, regardless of how the motion of agents should construct such distribution, is given in the seminal paper of Koopman \cite{koopman1957theory}. Assume for now that the target distribution $p$ is stationary. Koopman showed that the optimal coverage distribution up to time $t$ is 
\begin{equation}
c^t_{opt}=\max\left[\ln p-\alpha^t, 0\right],
\label{eq:optimal_coverage}
\end{equation}
where $\alpha^t$ satisfies
\begin{equation*}
Nt = \int_{S} \max\left[\ln p-\alpha^t, 0\right] \diff x.
\end{equation*}
From \eqref{eq:total_coverage} and \eqref{eq:total_smooth_coverage} it is clear that the total coverage is a linear function in time. Therefore, the total amount of available coverage in a time interval of size $\Delta t$ is equal to $N\Delta t$. Koopman also showed that the optimal solution can be achieved by incremental planning, i.e., if after spending $Nt$ of coverage, the target has not been found but some extra amount of coverage, $N\Delta t$, becomes available then we can use the same procedure, using the posterior probability of targets, to compute the optimal distribution for the next stage of the search, and yet the probability of detection would have been the same if we had allocated $N(t+\Delta t)$ from the beginning.

Although the optimal coverage distribution in \eqref{eq:optimal_coverage} can be readily computed, it is not known how to design the search agent paths to achieve this optimal distribution while satisfying physical constraints like the continuity of the paths. In fact, given any special form of coverage distribution, such as \eqref{eq:coverage} or \eqref{eq:c_smoothing}, achieving the optimal distribution might be infeasible. On the other hand, direct maximization of \eqref{eq:detection_prob} over the space of agent paths leads to a nonconvex problem which is highly dependent on the initial guess, and therefore too costly for real-time applications. In order to make this problem feasible, we reformulate the optimal search problem as follows. 

Assume that during the time interval $[0,t]$, $N$ agents have searched the area but the search target has not been detected yet. The current coverage distribution is given by $c^t$, which is not necessarily the optimal coverage. Now for the next stage of the search over the horizon $[t,t+\Delta t]$, there are $N_1$ agents available and the amount of available coverage is $ N_1 \Delta t$. If we assume that during the time period $\Delta t$ the drift of distribution $p^t$ is negligible when compared with dynamics of the search, the optimal coverage at time $t+\Delta t$ according to Koopman theory of search is given by \eqref{eq:optimal_coverage}, but this time $\alpha$ satisfies 
\begin{equation}
\int_{S} c^t_\sigma \diff x + N_1 \Delta t = \int_{S} \max \left[\ln p^{t} - \alpha^{t}, 0\right] \diff x.
\label{eq:alpha}
\end{equation}
Given that achieving this optimal coverage may be infeasible, \textit{we only strive to minimize the mismatch between the current distribution and the optimal distribution locally in time}. Therefore, the mDSMC algorithm uses a mismatch distribution which follows naturally from the Koopman search theory:
\begin{equation}
s_\sigma^t
= \max \left[ \ln p^t - \alpha^t - c_\sigma^t, 0\right].
\label{eq:source_modified}
\end{equation}
and we seek the next direction of motion for each agent that minimizes the scalar quantity
\begin{equation}
\Phi^t = \lVert s_\sigma^t \rVert
\label{eq:norm_modified}
\end{equation}

On the other hand, the DSMC algorithm mismatch distribution used in \cite{mathew2010uniform} has a much simpler form:
\begin{equation}
s^t
=  p^t - \frac{c^t}{Nt},
\label{eq:source}
\end{equation}
and we seek to minimize $\lVert s^t \rVert$.

There are some important advantages of using \eqref{eq:norm_modified} instead of  \eqref{eq:source}. The mDSMC mismatch excludes the “over-searched” areas i.e. when $c_\sigma^t >\ln p^t - \alpha^t$. which improves search results significantly and removes some instabilities which are often present with DSMC algorithm.

\section{Spectral Multi-scale Coverage path planning}

The key feature of spectral multi-scale algorithms is that they use a special class of function norms, called the Sobolev-space norm of negative index, to quantify the mismatch of current and optimal coverage. This norm for the mismatch distribution defined above is
\begin{equation}
\Phi^t = \sum_{k \in Z^2} \Lambda_k s^t_k,
\label{eq:phiseries}
\end{equation}
where $s^t_k$ is the coefficients of two-dimensional Fourier expansion of $s^t$, given by
\begin{equation*}
s^t_k = \int_{A}f_k s^t \diff x
\end{equation*}
where $f_k (x)$ are the Fourier basis functions, and $k$ is the wave-number vector. 
The norm coefficients in \eqref{eq:phiseries} are given by
\begin{equation*}
\Lambda_k = \left(1 + \lVert k \rVert^2 \right)^\beta,
\end{equation*}
where $\beta$ is the index of the norm and takes a negative value. This type of coefficient makes the contribution of components with smaller spatial scale (represented by larger $k$) proportionally smaller. Therefore, any algorithm that minimizes the mismatch function $\Phi^t$ automatically puts more weight on the large-scale spatial features of the target distribution. For the mDSMC algorithm $\beta=-\frac{1}{2}$ is used while for the DSMC $\beta\leq-\frac{3}{2}$ (in all DSMC calculations we used $\beta=-3/2$).
Lower exponent $\beta$ in mDSMC algorithm reduces smoothing of the mismatch function which becomes important later in the search when smaller spatial scale features are explored. 

DSMC and mDSMC algorithm minimize the mismatch distribution by implementing the instantaneous corrections to the agent paths that result in the fastest descend in the value of $\Phi$ at any time. Our modified algorithm can be extended to agents with second-order dynamics similar to \cite{mathew2010uniform} but here we assume that the agents have first-order dynamics with constant velocity. As shown in \cite{mathew2010uniform}, we first need to define the potential field
\begin{equation*}
u^t(z^t_i) = \sum_{k \in Z^2} \Lambda_k s^t_k f_k(x).
\end{equation*}
Then the velocity vector for each agent is given by 
\begin{equation*}
v^t_i = \frac{\nabla u^t(z^t_i)}{\lVert \nabla u^t(z^t_i) \rVert} v_{mag,i}
\end{equation*}
where $v_{mag,i}$  is the constant velocity magnitude of the $i$-th agent. 
Finally, the movement of the search agents is governed by the first order motion law:
$\frac{\diff z^t_i}{\diff t} = v^t_i$.

To summarize, the mDSMC in contrast to DSMC uses a new type of objective function, described in \eqref{eq:source_modified} and minimizes the Sobolev norm with index $\beta=-1/2$, which is better suited to search problems.
In this specific application of mDSMC algorithm for the RBF, we used a compact support approximation of normal distribution centered at the agent location with a standard deviation of $\sigma = 3 \text{km}$ and 1-hour window to determine the amount of available coverage in \eqref{eq:alpha}.

\section{Search for MH370 (Results)}

The search and rescue operation for MH370 began on March 9, 2014, one day after the loss of communications between the plane and ground stations. The surface search lasted for 50 days and included areas near the Malay Peninsula and the Southern Indian Ocean. The estimated splash location was constantly updated due to a stream of incoming information and the drifted image of the most probable splash area was drifted under the ocean model to obtain the target area for each stage of the search.  We consider three probable splash areas: 
\begin{itemize}
	\item Area A: the drifted image of this area was searched from March 28 to April 1,
	\item Area B: the drifted image of this area was searched on April 2 and 3,
	\item Area P: this area, in its totality, was not searched in the surface search. However, after the conclusion of the surface search, area P was recognized as the most probable splash area and marked as the priority location for underwater search \cite{bureau2017operational}.
\end{itemize}

In the actual search for MH370, the probable splash areas (including A and B above) were drifted using the knowledge of the ocean currents and wind effects, over the time lapse between the splash and commencement of the search. A precise description of the used drift model is not disclosed; however, the coordinates of the splash areas and some of the drifted images (including images of A and B) are publicly available. In our computational model, we use the surface velocity data obtained from the HYbrid Coordinate Ocean Model (HYCOM) to drift the splash areas. Moreover, we assume that the initial target distribution is uniform over each area. 

We simulate the following scenarios:
\begin{enumerate}
	\item \textit{Lawnmower scenario 1}: the reported search areas are searched using the lawnmower algorithm representing the strategy that occurred in the actual MH370 search.
	\item \textit{Lawnmower scenario 2}: the splash areas are drifted using our drift model, and then convex-hull of target distribution is searched using the lawnmower algorithm.
	\item \textit{DSMC}: the splash areas are drifted using our drift model, and then searched using the DSMC algorithm.
	\item \textit{mDSMC}: The splash areas are drifted using our drift model, and then searched using the mDSMC algorithm.
\end{enumerate}

Note that scenario 1 is a replication of the actual search. Contrasting scenario 2 and 3 provides a comparison of the search algorithms independent of the drift model, while contrasting scenario 1 and 3 distinguishes the combined effect of our drift model and search algorithm.  
Comparing scenario 3 and 4 shows the advantages gained by modifying DSMC.

\subsection{Problem Definition of Search for MH370}

The search parameters used in this study are based on two main sources: The coordinates of the probable splash areas are extracted from a series of reports on MH370 search operation by the Australian Transport Safety Bureau (ATSB) \cite{bureau2017operational}, and the information on search agents and the searched areas are collected from media releases and daily briefings by the Joint Agency Coordination Centre (JACC) for MH370 search. The target areas of search by JACC were determined by drifting the splash areas using models of the surface currents complemented with the real-time wind and wave data, however, the details of the drifting model are not publicly available. The search areas specified by JACC are used in scenario 1. In scenario 2, we have used the splash areas reported by ATSB and drifted them using our model. The search area for lawnmower was then computed as the smallest convex bounding polygon around the area with a nonzero probability of targets.

Up to 21 aircraft and 19 ships were deployed in the search operation. In our study, we have only used the aircrafts as search agents, due to the lack of data on the technical specifications of the ships.  The number and type of the deployed aircraft vary with the search day (see Table~\ref{tbl:search_parameters}). but we have chosen the speed of every agent to be 380 km/h which represents the typical loiter speed of the military aircraft involved in the search. We also assume that the scanning of the search areas on each day started at 2:00 pm and ended at 5:00 pm (UTC).

\begin{table}
	\caption{The parameters of MH370 search operation.}
	\label{tbl:search_parameters}
	\def\arraystretch{1.3}
	\begin{tabular*}{\columnwidth}{@{\extracolsep{\stretch{1}}}*{3}{c}@{}}
		\hline 
		Search Area 			& Seach Period 			& Number of search aircraft \\ 
		\hline 
		Area A      			& March 28 - April 1 	& 10, 8, 10, 10, 11 \\ \\
		Area B      			& April 2 - 3 			& 9, 9 \\ \\
		\multirow{3}{*}{Area P}	& March 8-17,			& 	\multirow{3}{*}{10, 10, 10, 10, 10, 10, 10, 10, 10, 10}	\\ 
		& March 13-22, &  \\ 
		& March 18-27 &  \\ 
		\hline 
	\end{tabular*} 
\end{table}

A simplified stochastic model is used to emulate the target detection by the search observers aboard the aircraft. In this model, if the target remains within the 1.5-km radius of an agent for time $t$, then it will be detected with the probability $P=1-\exp(-t/T)$, where $T=2$ seconds is the expected detection time. 
A machine vision system for a low-cost fixed-wing UAV has been investigated in \cite{leira2015automatic} where thermal imaging camera and onboard processing unit are used for performing real-time detection, classification, and tracking of objects floating on the ocean surface. Bearing in mind that mainly an eye vision is used in MH370 search, the authors believe the expected detection time of 2 seconds used in simulations is a reasonable estimate.

\subsection{Drift model}
Our drift model is based on the surface velocity data computed by the US Navy, using the HYbrid Coordinate Ocean Model (HYCOM). The data consists of 3-hourly longitudinal and latitudinal velocity components with a spatial resolution of 1/25 degrees in each direction. The drift model is used to compute the flow map $\mathcal{T}$ i.e. the paths of the search targets as well as the trace of search agent paths in \eqref{eq:coverage}. The location of the targets in the splash area is initialized at 0:30 UTC on March 8 using the Halton sequence and then updated through a 4th-order adaptive Runge-Kutta method with linear interpolations of the surface velocity in time and space. 
The evolution of the target distribution is simulated using a semi-Lagrangian method: we uniformly sampled the probable splash areas with a high number of tracers ($10^5$ tracers per area) and drifted those tracers using the above drift model. By computing the local density of tracers, we assemble the distribution at future times. This method does not require a numerical mesh and eliminates the artificial diffusion associated with Eulerian schemes for solving \eqref{eq:drift}.

\begin{figure}[!htb]
	\centering
	\includegraphics[width=0.5\linewidth,keepaspectratio=true]{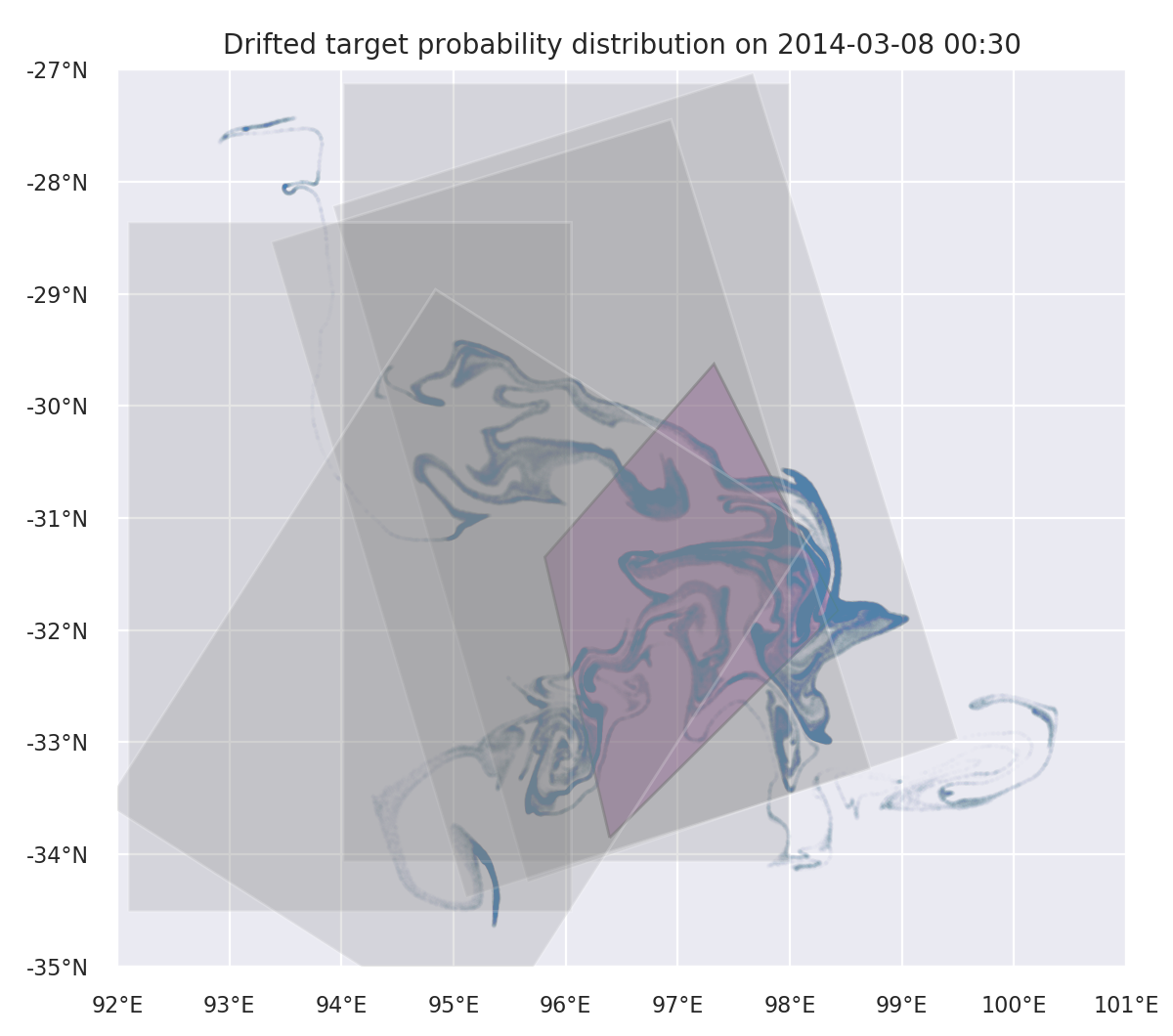}
	\caption{The estimated splash area A (purple diamond) and the drifted probability distribution of search targets on the first day of actual search (March 28).  Gray polygons represent regions of actual search which is simulated in Lawnmower scenario 1 which was conducted from March 28 to April 1.}
	\label{fig:area_a_splash}
\end{figure}

\begin{figure}[!htb]
	\centering
	\includegraphics[width=0.5\linewidth,keepaspectratio=true]{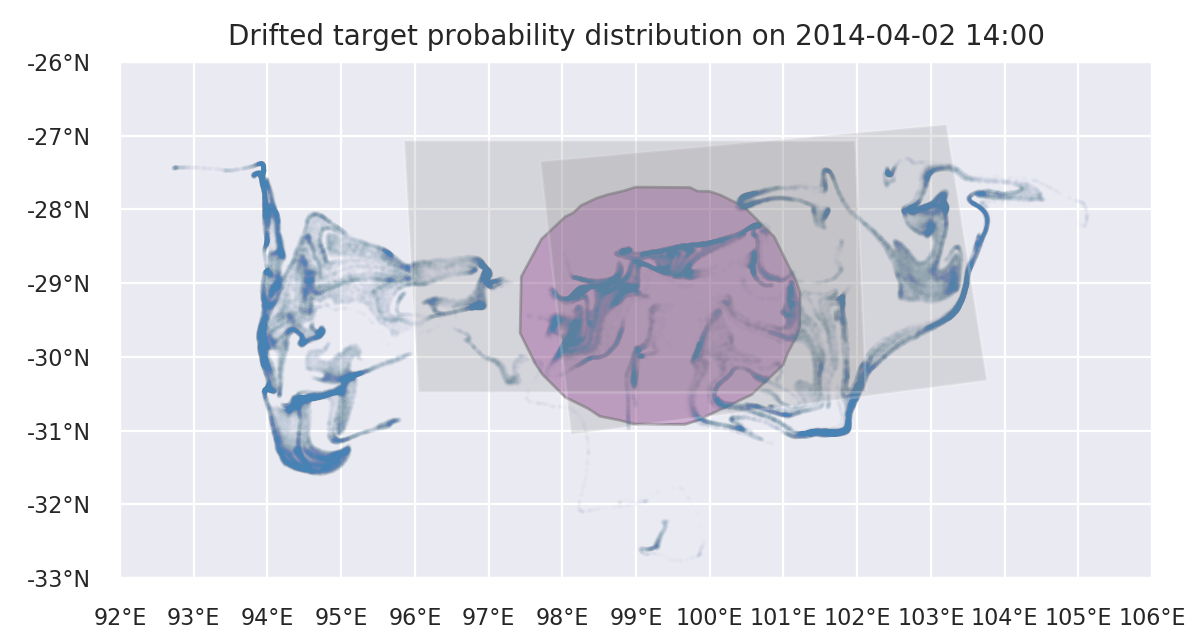}
	\caption{The estimated splash area B (purple circle) and the drifted probability distribution of search targets on the first day of actual search (April 2). Gray polygons represent regions of actual search which is simulated in Lawnmower scenario 1 which happened on April 2 and April 3.}
	\label{fig:area_b_splash}
\end{figure}

Figure~\ref{fig:area_a_splash} shows the drift of the area A. During the 19-day lapse between the estimated splash time and the start of the search in this area, the probability distribution of the search targets undergoes substantial change, and a potential search target moves several hundred kilometers away from its initial position (purple diamond). The reported areas of actual search are mainly covering the drifted target probability, but it is evident that a lot of search effort is spent on regions where target probability is zero. For area B, a similar visualization is shown in Figure~\ref{fig:area_b_splash}. Seemingly, the conducted search regions are well placed, but a large part of the target probability in the west is not covered at all.  
The stretching and folding of the support of the target distribution in Figures \ref{fig:area_a_splash} and \ref{fig:area_b_splash} is typical of unsteady geophysical flows with strong mixing behavior \cite{mezic2010new}. 
Relative dispersion properties in the ocean surface layer further complicate drifting estimation \cite{corrado2017general}.
Using a lawn-mower search algorithm, which is suited to regular shapes, would be inefficient in such situations. In contrast, the paths of the mDSMC would naturally adapt to the shape of the drifted area. 

\subsection{The search}

\begin{figure}[!htb]
	\centering
	\includegraphics[width=\linewidth,keepaspectratio=true]{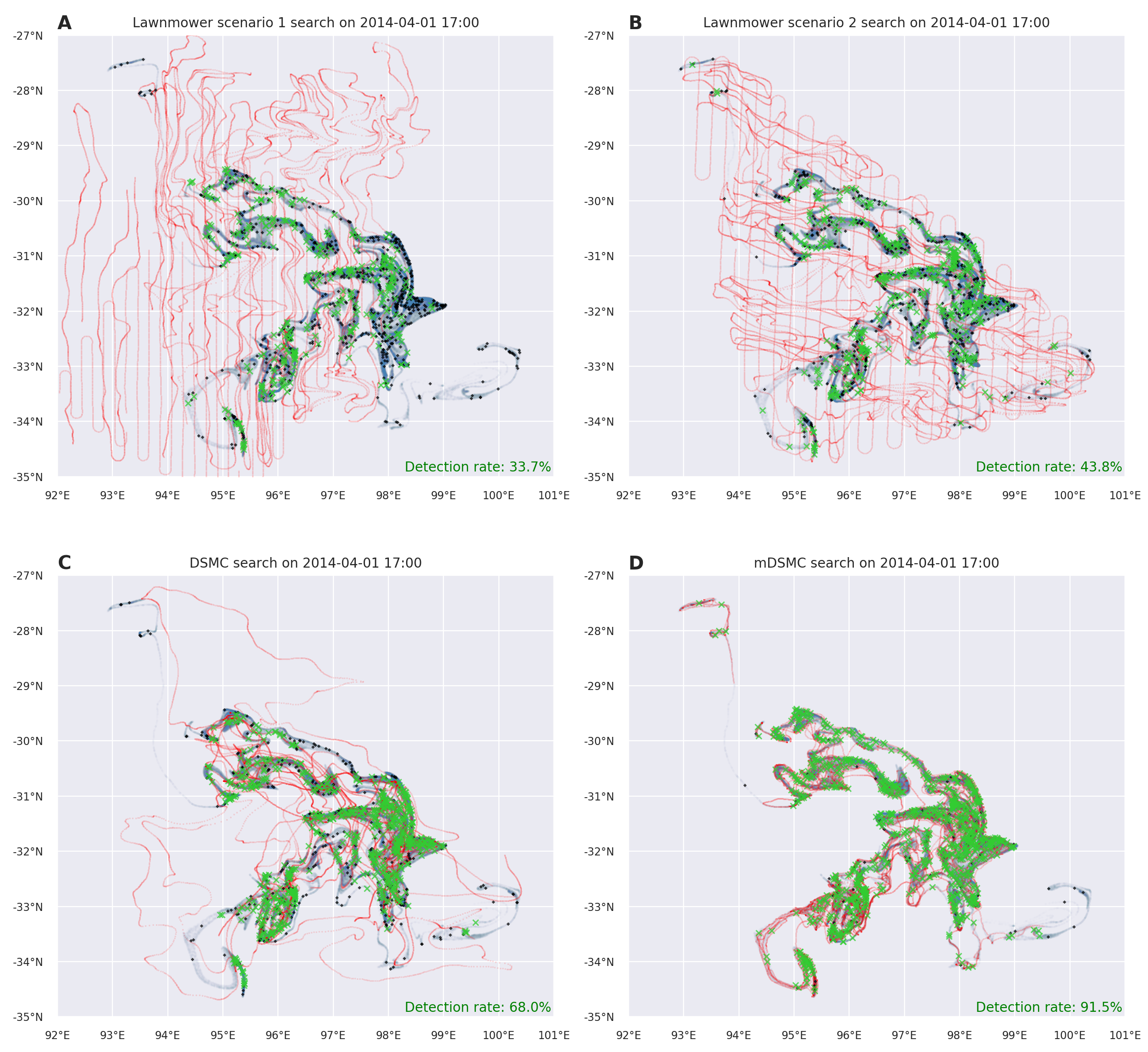}
	\caption{Path of agents in different search scenarios for MH370: (A) The estimated splash area A at March 9 (diamond) and the drifted probability distribution of search targets on the first day of actual search (March 28). (B) The paths of search agents in the mDSMC search on the drifted image (scenario 3) (C) The lawnmower path of search agents in the reported search (scenario 1). (D). The path of search agents in lawnmower search on the drifted image (scenario 2). }
	\label{fig:area_a_search}
\end{figure}

The visualization of all four considered scenarios for area A search is shown in Figure~\ref{fig:area_a_search}. The lawnmower areas of actual MH370 search yield the trajectories shown in \ref{fig:area_a_search}.A, while \ref{fig:area_a_search}.B show trajectories of lawnmower strategy which relies on our drift model. Both sets of trajectories fails to achieve acceptable results in terms of detection rate (33.7\% and 43.8\%, respectively) which is indicated with green crosses representing detected targets.

The DSMC and mDSMC search strategies (Figures \ref{fig:area_a_search}.C and \ref{fig:area_a_search}.D, respectively) are guided by the drifted target probability distributions and, as result of that, accomplished trajectories are passing through the regions of high probability and consequentially produce better target detection rate. 
The improvements made in mDSMC, in contrast to DSMC, are recognizable in a better local search and in the elimination of wide and unnecessary paths. It is rather interesting to observe how mDSMC covers thin ribbons of the target probability distribution, especially the one in the north of the search domain.

\begin{figure}[!htb]
	\centering
	\includegraphics[width=\linewidth,keepaspectratio=true]{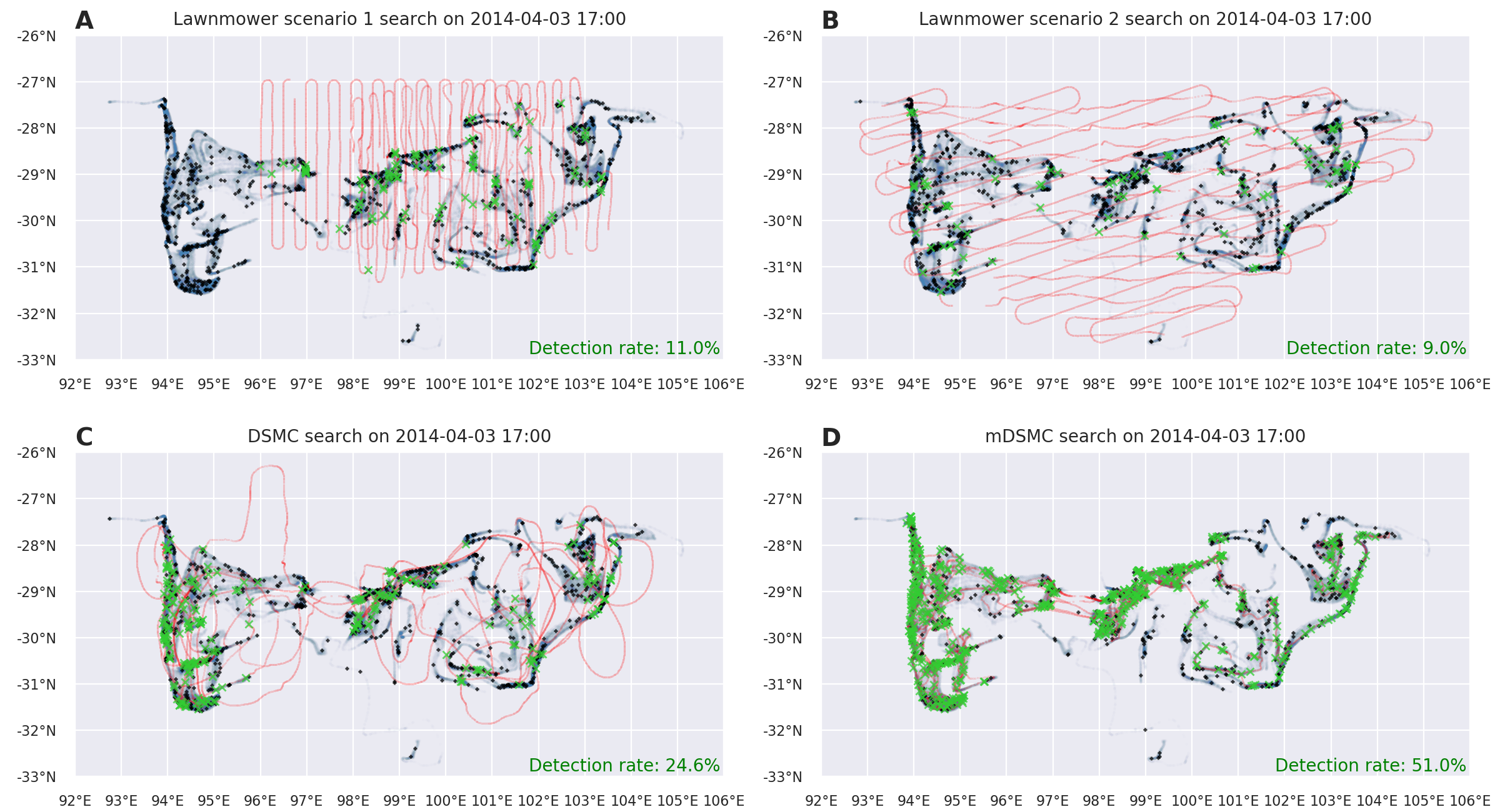}
	\caption{Path of agents in different search scenarios for MH370: (A) The estimated splash area B at March 9 (diamond) and the drifted probability distribution of search targets on the first day of actual search (April 2). (B) The paths of search agents in the mDSMC search on the drifted image (scenario 3) (C) The lawnmower path of search agents in the reported search (scenario 1). (D). The path of search agents in lawnmower search on the drifted image (scenario 2).}
	\label{fig:area_b_search}
\end{figure}

Figure~\ref{fig:area_b_search} displays the results of search simulations for area B.
In lawnmower strategies, an unexpected anomaly can be observed: the lawnmower scenario 2, shown in Figure~\ref{fig:area_b_search}.B, although taking into account probability distribution drifted with the same drift model as targets, is apparently not better than the "naive" lawnmower scenario 1 ( Figure~\ref{fig:area_b_search}.A) if the number of detected targets is compared. Even though the lawnmower scenario 1 strategy missed a big portion (approx. 1/4) of targets on the west part of the domain, due to the density of lawnmower trajectories in remaining regions, it outperforms the scenario 2 which sparsely covers whole target distribution. Similar to area A results, both DSMC and mDSMC (Figures \ref{fig:area_b_search}.C and \ref{fig:area_b_search}.D, respectively) achieved better results than lawnmower strategies. Due to the highly indented and complex shape of the target probability distribution, DSMC's tendency to global search prevents more accurate search.

To evaluate the success of different search scenarios, we have performed the following experiment: For each case, we have seeded the splash area on the day of the splash with 1000 randomly positioned targets and drifted the targets with the surface flow. On each day of the search, the search algorithms are executed while the targets are drifting during the search, and the success rate of each algorithm is computed as the fraction of targets detected by the agents. The detection is modeled as a stochastic process with the law of exponential saturation. The described procedure of a search simulation is repeated 100 times for each scenario using randomly varying initial positions for the targets and agents, and the ensemble average of the success rates are reported here.

\begin{figure}[!htb]
	\centering
	\includegraphics[width=\linewidth,keepaspectratio=true]{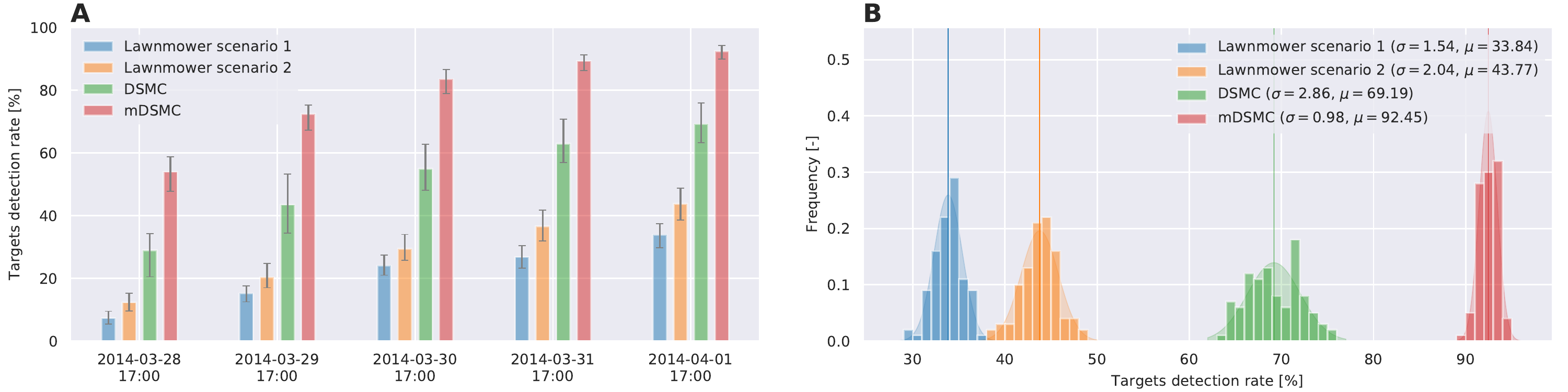}
	\caption{Success rate of search scenarios for MH370:
		(A) The progress of the search scenarios vs time for drifting area A. The target detection only moves upward in the time interval of 14:00-17:00 on each day. This constraint is posed by the fact that the search area is about 700 km away from the operation base.
		(B) Histogram of the success rate for search scenarios in 100 realizations of each scenario with random initial conditions. }
	\label{fig:area_a_comparison}
\end{figure}

\begin{figure}[!htb]
	\centering
	\includegraphics[width=\linewidth,keepaspectratio=true]{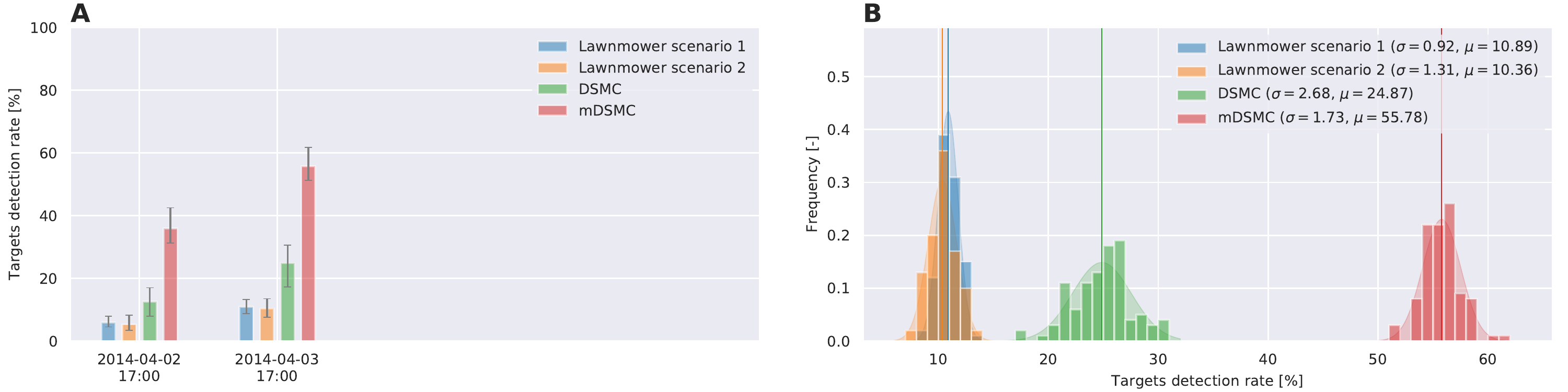}
	\caption{Success rate of search scenarios for MH370: (A) The progress of the search scenarios vs time for drifting area B. The target detection only moves upward in the time interval of 14:00-17:00 on each day. This constraint is posed by the fact that the search area is about 2700 km away from the operation base.  (B) Histogram of the success rate for search scenarios in 100 realizations of each scenario with random initial conditions. }
	\label{fig:area_b_comparison}
\end{figure}

Figures~\ref{fig:area_a_comparison} and \ref{fig:area_b_comparison} summarize the results of the above experiment for area A and area B, respectively. The execution of the search algorithms is limited to a 3-hour period each day to account for the travel time of the agents between the closest airbase and the search area. The performance of DSMC and mDSMC algorithms is remarkably different from the lawnmower technique: the lawnmower scenarios have uniform daily progress which is expected since the lawnmower covered area progresses linearly with time. 
The DSMC generates trajectories that roughly follow the target distribution - addressing the broad non-local search which is improved in mDSMC.
In the case where mDSMC algorithm is used, however, the search agents perform large sweeps of the area first, which results in capturing about \%50 of the targets on the first day. In the subsequent two days, the agents perform increasingly localized search to find the remaining targets. The results for both areas A and B show that mDSMC algorithm would have significantly improved the chances of finding floating debris from the MH370 airplane.
The dynamics and comparison of the search for area A and B can be found in Video 1 and Video 2 included in Supplementary.

\section{Effect of delayed search on the success rate}

An unfortunate circumstance of the MH370 search was that a reliable estimate of the splash location was not known immediately after the plane’s disappearance. The day after the plane went missing, the search started by focusing on the vicinity of the Malay Peninsula. In the meantime, the analysis of satellite data indicated that the plane had continued to fly for six hours after the last communication, and two preliminary corridors, one pointing north toward Kazakhstan and one to the south toward the Indian Ocean, were speculated as to the probable flight paths. By March 24, a consensus was formed around the hypothesis that the plane had crashed on the ocean surface in the southern Indian Ocean. 

\begin{figure}[!htb]
	\centering
	\includegraphics[width=\linewidth,keepaspectratio=true]{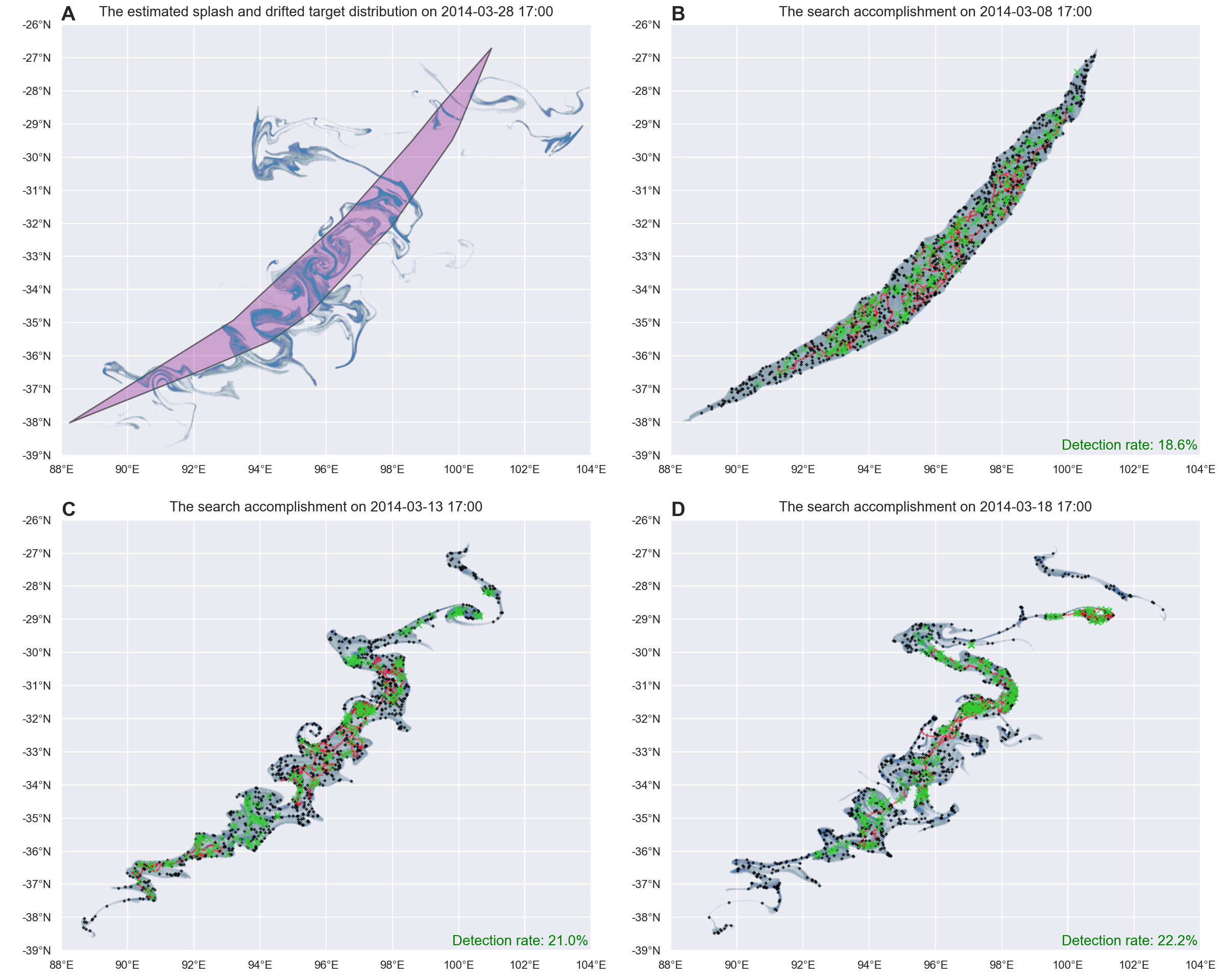}
	\caption{(A) The estimated splash area P (purple region) and the drifted area P for March 28. (B, C, D) show drifted target distribution, the agent paths and target detection at the end of the first search days as a result of search simulation starting on the day of the splash, 5 and 10 days after, respectively.}
	\label{fig:area_p_search}
\end{figure}

An important question that arose in these circumstances is how the dispersion of the debris until the start of the search would have affected the success of the search performance.
To answer this question, we have simulated the mDSMC search on the area P starting 0, 5 and 10 days after the day of the splash. 
Figure~\ref{fig:area_p_search} illustrates the dynamics of ocean surface mixing and the complexity of drifted target distribution. The initially estimated splash area and its drifted image after 20 days is shown in (Fig.~\ref{fig:area_p_search}.A). The drifted target probability distribution and the accomplishment of the search started on the day of the splash, 5 and 10 days later are shown in figures \ref{fig:area_p_search}.B, \ref{fig:area_p_search}.C and \ref{fig:area_p_search}.D, respectively. As time goes by, the target distribution is accumulated in certain parts of the search domain and, since it is aware of the target distribution, mDSMC control method guided agents towards parts with higher target probability. At the start of the search (Figure~\ref{fig:area_p_search} show searches after a day of the search) trajectories are more accumulated as the later search start time.

\begin{figure}[!htb]
	\centering
	\includegraphics[width=\linewidth,keepaspectratio=true]{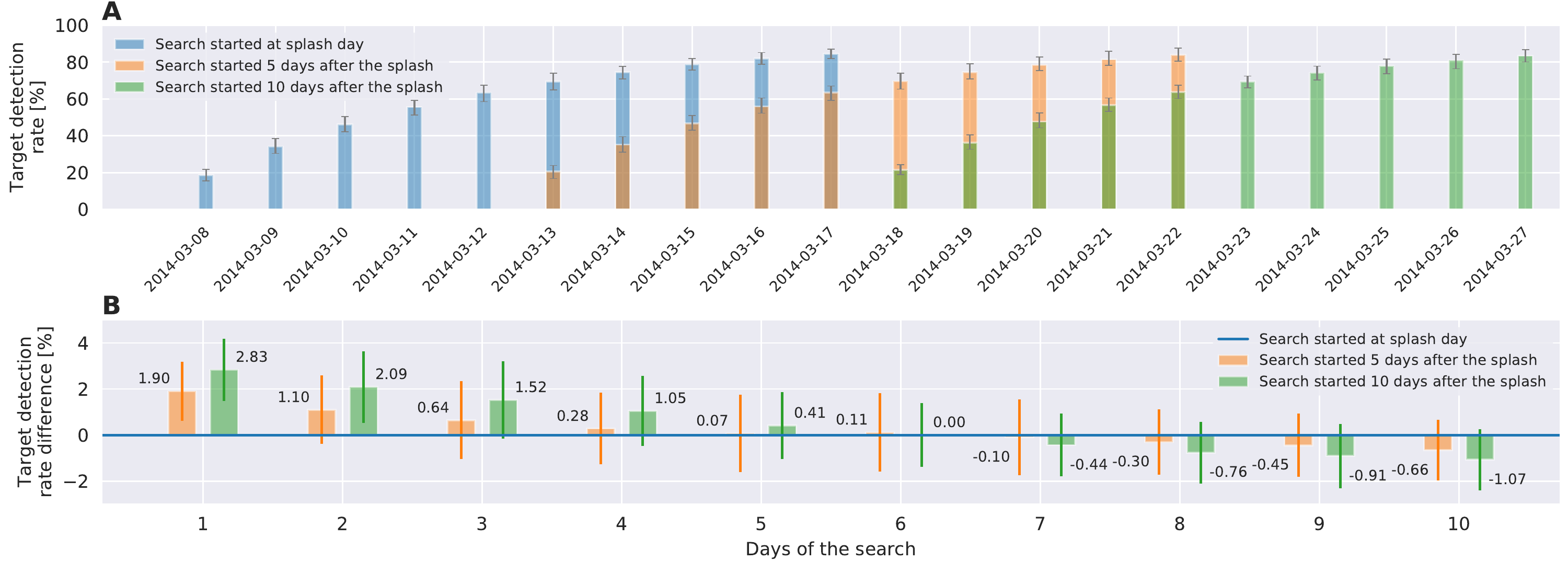}
	\caption{Effect of delayed search on the success rate of mDSMC algorithm: (A) The progress of the mDSMC search algorithm on area P. (B) Variations in the success rate of daily search due to a 5-day or 10-day delay in the commencement of search.}
	\label{fig:area_p_comparison}
\end{figure}

Figure~\ref{fig:area_p_comparison}.A shows the daily progress of the search algorithm starting with 0-day (immediately after the splash), 5-day and 10-day delay, while Figure~\ref{fig:area_p_comparison}.B depicts the variations thereupon for the search operations starting 5 or 10 days later. The results, surprisingly, indicate that the chances of finding debris in a 5-day search are greater if the search is started later.
The consistency of the detection rate over 100 runs indicates the robustness of the mDSMC method.

\begin{figure}[!htb]
	\centering
	\includegraphics[width=\linewidth,keepaspectratio=true]{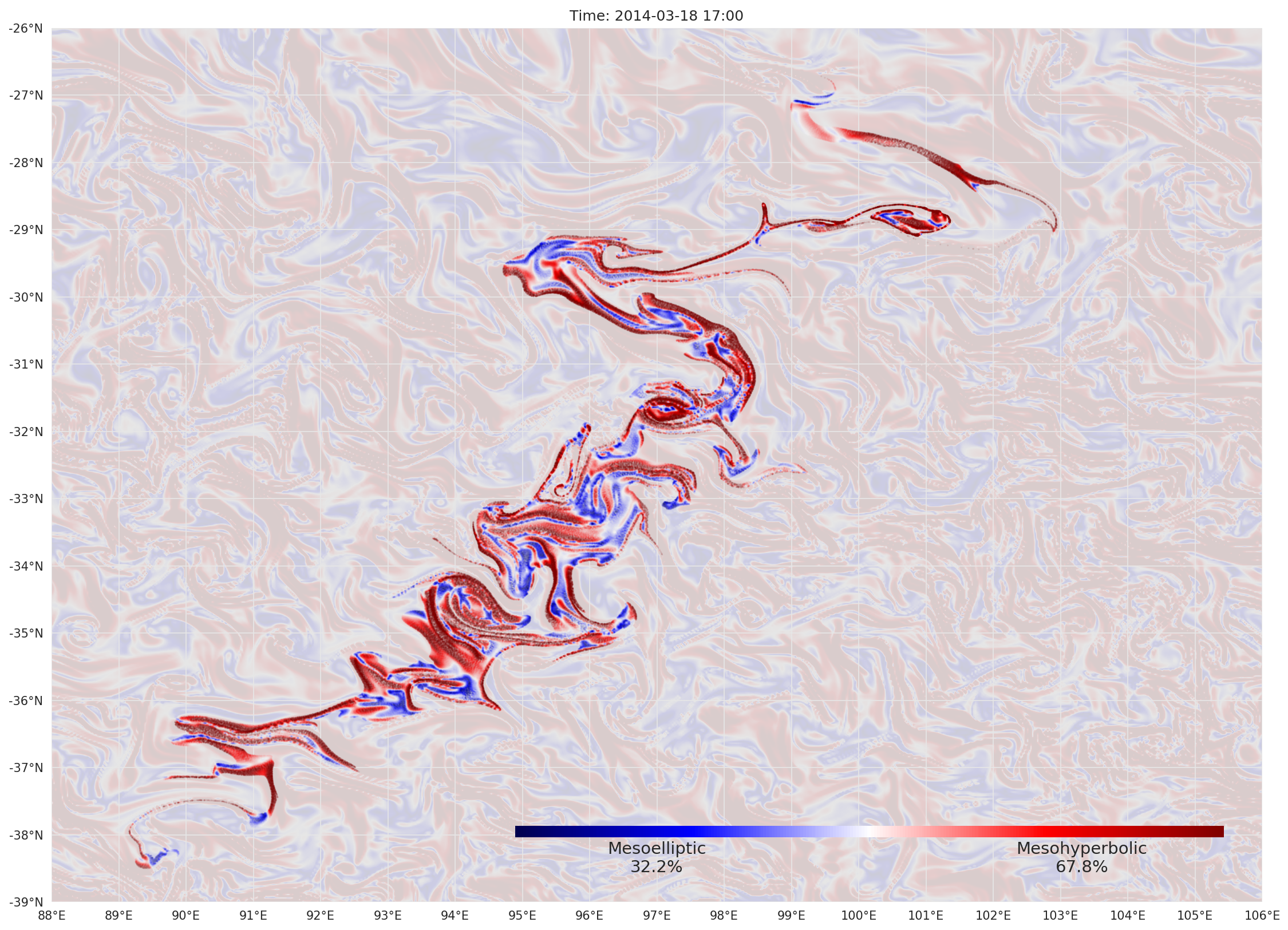}
	\caption{Finite-time mixing analysis of the search area: The colors depict the finite-time Lagrangian behavior in the drifted area A. Dominant rotational behavior is marked in blue and stretching behavior in red. The search targets accumulate around red ridges which reduces the effective size of the search area.}
	\label{fig:area_p_hypergraph_2zones}
\end{figure}

The above observation can be explained by inspecting the surface mixing in the search area. To do so, we employ the so-called hypergraph analysis introduced in \cite{mezic2010new}. The hypergraphs provide qualitative maps of finite-time motion for small objects carried by unsteady flows. Based on determinant $\det \left(\nabla \mathcal{T}^{t_1, t_2} / (t_2, t_1)\right)$ the flow domain is partitioned into two qualitative classes of behavior: mesoelliptic regions where the motion is dominated by rotation, and mesohyperbolic regions where motion is dominated by stretching in one direction and contraction in the other. The hypergraph of the area P for mixing over the 12-day interval starting at the splash day is shown in Figure~\ref{fig:area_p_hypergraph_2zones}, with hyperbolic behavior, marked in red and elliptic in blue. The key observation in \cite{mezic2010new} was that hyperbolic regions reveal the convergence zones for floating objects on the surface. As shown in the figure, the majority of the search targets would accumulate around these regions over the first few days after the splash. As a result, the delayed search operations would initially encounter a higher concentration of targets in those areas and therefore yield higher success rates at the beginning. 

\section{Conclusion}

The mDSMC algorithm is shown to be a viable and efficient algorithm for real-time planning of search operations in a dynamic and large environment such as ocean surface. This algorithm is based on a feasible formulation of the classic optimal search theory and is capable of addressing complex geometries with large uncertainty that arise in complex and unsteady environmental dynamics. Important features of mDSMC are the awareness of search agents of the large-scale spatial structure of the target area, the continuous balance of search load between agents and the ability to concentrate the search on most feasible areas. mDSMC brings an improvement over DSMC in local search, where high- detailed structures of the probability distribution are explored. The application of this algorithm to the MH370 case showed that these features lead to a several-fold improvement in the success rate of the search compared to the conventional lawnmower method. This suggests that using this algorithm could lead to a critical improvement in the survivor rescue operations and/or finding critical evidence in similar incidents. The mDSMC algorithm can also be employed in other forms of operations that involve mobile agents (planes, drones, land patrols, etc.) and moving targets (individuals, animals, debris, black box, etc.). Examples of such operations include geographical and zoological surveys, underwater sonar mapping and search, the rescue of/from wildlife, contamination removal and spill containment in natural reserves and oceans.  

A surprising conclusion of our study was that delaying the search operation in dynamic environments can lead to higher success rates in finding the debris in the initial stages of the search. In the case of MH370, this is caused by the convergence zones within the initial uncertainty distribution. Over the first few days, these zones attract and accumulate the target probability density on the surface. A search operation at the right time can exploit this accumulation of density and lead to a higher success rate given the same amount of coverage. The convergence property is typical of a dynamic environment that shows regular patterns for the accumulation of targets (e.g. shrinkage of area, or existence of attracting regions), and it can be exploited in other problems to increase the efficiency of search and sampling in terms of search time and the number of agents.

\bibliography{bibliography}

\section*{Acknowledgements}
The authors are grateful to Pat Hogan and Ole Martin Smedstad for providing the Global Ocean Prediction data for the drift model. Funding: B.C.’s research is fully supported by the Croatian Science Foundation under the project IP-2019-04-1239. S. I.’s contribution is supported by the University of Rijeka under the project number 17.10.2.1.04. I.M's research is supported by ARO-MURI grant W911NF-17-1-0306. 

\section*{Author contributions}
I.M. conceived the idea for the article. S.L. and P.C. verified the data sets and prototyped the search algorithm. S.I. and B.C. proposed the modifications to the search algorithms and simulated the experiments. H.A. collected the search parameters, formulated the search problem, and wrote the initial manuscript. 

\section*{Additional information}

All authors were involved in analyzing the results and final edition of the manuscript. Competing interests: The authors declare that they have no competing interests. Data and materials availability: All parameters for reproducing the study are presented in the paper. Additional data related to the paper may be requested from the authors.

\section*{Competing interests}
The authors declare no competing interests.

\section*{Video Captions}

\begin{itemize}
	\item[] Video 1: Simulation of area A search
	\item[] Video 2: Simulation of area B search
\end{itemize}

\end{document}